\documentclass [reqno,10pt,oneside]{amsart}
\usepackage{amsthm}
\usepackage{amscd}
\usepackage{amsmath,amssymb}
\usepackage{algorithm}
\usepackage[noend]{algorithmic}

\newcommand{\singular}{{\sc Singular }}
\newcommand{\bertini}{{\sc Bertini }}
\usepackage{tikz}
\usepackage{amsthm}
\usetikzlibrary{decorations.pathreplacing}
\usetikzlibrary{arrows}
\newtheorem{rem}{Remark}
\newtheorem{eg}{Example}
\newtheorem{theo}{Theorem}

\newtheorem{defi}{Definition}

\numberwithin{equation}{section}
\numberwithin{theo}{section}
\numberwithin{prop}{section}
\numberwithin{lem}{section}
\numberwithin{cor}{section}
\numberwithin{defi}{section}
\numberwithin{rem}{section}
\numberwithin{eg}{section}
\numberwithin{alg}{section}
\begin{document}

\title{Numerical Decomposition of Affine Algebraic Varieties}
\author{Shawki AL-Rashed \;\;\;\;\;\;\; Gerhard Pfister\\ \;\;rashed@mathematik.uni-kl.de\;\;\;\; pfister@mathematik.uni-kl.de\\\\\\
Department of Mathematics,\\ University of Kaiserslautern,\\Erwin-Schr\"odinger-Str.,\\
    67663 Kaiserslautern,\\ Germany}
\date{\today}
\maketitle
\thanks{}

\begin{abstract}
\; An irreducible algebraic decomposition $\cup_{i=0}^{d}X_i=\cup_{i=0}^{d}(\cup_{j=1}^{d_i}X_{ij})$ of an affine algebraic variety X can be  represented as an union of finite disjoint sets $\cup_{i=0}^{d}W_i=\cup_{i=0}^{d}(\cup_{j=1}^{d_i}W_{ij})$ called numerical irreducible decomposition (cf. \cite{Sommese1},\cite{Sommese2},\cite{Sommes2},\cite{Sommes3},\cite{Sommese3},\cite{Sommese5},\cite{Sommese6},\cite{Sommes4}). $W_i$ corresponds to a pure i-dimensional $X_i$, and $W_{ij}$ presents an i-dimensional irreducible component $X_{ij}$.
Modifying this concepts by using partially Gr\"{o}bner bases, local dimension, and the "Zero Sum Relation" we present in this paper an implementation in SINGULAR to compute the numerical irreducible decomposition. We will give some examples and timings, which show that the modified algorithms are more efficient if the number of variables is not too large. For a large number of variables \bertini is more efficient. Note that each step of the numerical decomposition is parallelizable. For our comparisons we did not use the parallel version of \bertini .
\end{abstract}
\textbf{keyword}: Witness point set, Homotopy function, Gr\"{o}bner basis, Local dimension, Monodromy action, Zero Sum Relation.
\section{Introduction}
 \; Given a system of n polynomials in $ \mathbb{C}^N$,
 \begin{equation*}
  f(x_1,...,x_N):=\left(
                   \begin{array}{c}
                     f_1(x_1,...,x_N) \\
                     . \\
                     . \\
                     f_n(x_1,...,x_N) \\
                   \end{array}
                 \right).
 \end{equation*}
Let X=V(f) be the algebraic variety defined by the system above. X has an unique algebraic decomposition into d pure i-dimensional components $X_i$,  $X=\cup_{i=0}^{d}X_i$. Where $X_i=\cup_{j_i}^{d_i}X_{ij}$ is the union of $d_i$ i-dimensional irreducible components, $d_0, d_1,...,d_d$ positive integers.\\
   The numerical irreducible decomposition (cf. \cite{Sommese2},\cite{Sommes2},\cite{Sommes3},\cite{Sommese3},\cite{Sommese6}) is given as the union $W=\cup_{i=0}^{d}W_{i}=\cup_{i=0}^{d}(\cup_{j=1}^{d_i}W_{ij})$. The $W_i$ are called i-Witness point sets and are given as an intersection of the pure i-dimensional component $X_i$ of X with a generic linear space L in $\mathbb{C}^N$ of dimension N-i, the $W_{ij}$ are called the irreducible witness point sets presenting the irreducible components $X_{ij}$ of dimension i with the following properties:
  \begin{itemize}
   \item $W_{ij}$ consists of a finite number of points contained in $X_{ij}$.
   \item $\sharp(W_{ij})=$deg($X_{ij}$) for $i\neq0$.
   \item $W_{ij}\cap W_{il}=\emptyset$ for $j\neq l$.
  \end{itemize}
 Computing of the numerical irreducible decomposition uses numerical polynomial homotopy continuation methods(cf. \cite{Verschelde1},\cite{Verschelde2}). This requires that the number n of a given polynomial system has to be equal to the number N of the variables. Therefore we reduce the polynomial system which defines X to a square system of N polynomials in N variables (cf. \cite{Sommese2},\cite{Sommes2},\cite{Sommese6}).
 Numerical irreducible decomposition (cf. \cite{Sommese2},\cite{Sommes2},\cite{Sommes3},\cite{Sommese3},\cite{Sommese6}) is proceeded in three steps:
 \begin{itemize}
  \item $1^{st}$ step reduces the polynomial system to a system of N polynomials in N variables and computes a finite set $\widehat{W}_i$ called witness point super set for a pure i-dimensional component $X_i$ for i=0,...,N-1. $\widehat{W}_i$ consists of points on $X_i$ and $J_i$ a set of points on components of larger dimension the so-called  Junk point set (cf. \cite{Sommese2},\cite{Sommese6}).
  \item $2^{nd}$ step removes the points of $J_i$ from $\widehat{W}_i$ to obtain a subset $W_i$ of the pure i-dimensional component $X_i$ (cf. \cite{Sommese6}).
  \item $3^{rd}$ step breakup $W_i$ into irreducible witness point sets representing the i-dimensional irreducible components of X using two algorithms. The first algorithm finds points on the same irreducible component in the Witness point set connected by path tracking techniques applying the idea of monodromy. The second algorithm computes a linear trace for each component, which certifies the decomposition (cf. \cite{Sommese1},\cite{Sommese5}).
 \end{itemize}

 In the second section we present a modified algorithm to compute $\widehat{W}_i$ using the cascade algorithm (cf. \cite{Sommese2},\cite{Sommes2},\cite{Sommese6}) and Gr\"{o}bner bases in the zero-dimensional case. Then we use the homotopy function (cf. \cite{Sommese6}), local dimension and Gr\"{o}bner bases in the zero-dimensional case to remove Junk points from $\widehat{W}_i$ to obtain the i-Witness point set $W_i$. This is explained in third section. In the fourth section we explain how to use the "Zero Sum Relation" and monodromy action on the algebraic variety to breakup $W_i$ into irreducible witness point sets. In the fifth section some examples are tested and timings are given on the basis of our \textsc{Singular} implementation.
\section{Witness Point Super Set}\label{WitSupSet}
 \begin{defi}
 Let Z=V(f) be an affine algebraic variety in $\mathbb{C}^{N}$ of dimension d, and X be a pure i-dimensional component of Z. Let $L_{i}$ be a generic linear space in $\mathbb{C}^{N}$ of dimension N-i.
 A set $\widehat{W}_{i}\subset\mathbb{C}^{N}$ is called i-Witness Point Super Set for X if it has the following properties:
  \begin{itemize}
     \item $\widehat{W}_{i}$ is a finite set of points.
     \item $ X\cap L_{i}\subset\widehat{W}_{i}\subset Z\cap L_{i}$ .
    \end{itemize}
\; The union $\widehat{W}$ of all i-Witness Point Super Sets is called a Witness Point Super Set for Z.
\end{defi}
\; In (cf. \cite{Sommese2},\cite{Sommes2},\cite{Sommese6}) the cascade algorithm is used to compute $\widehat{W}_i$. It starts with i=N-1 to compute the Witness point Super Sets $\widehat{W}_{i}$. It needs to define a start system G(x)=0 for the homotopy continuation method (cf. \cite{Verschelde1},\cite{Verschelde2}) and needs to know its solutions. We use a Gr\"{o}bner basis to compute the dimension d of Z, then use the cascade algorithm (cf. \cite{Sommese2},\cite{Sommes2},\cite{Sommese6}) which starts with i=d-1. We show below that we do not need to define a start system.
\begin{algorithm}
\caption{\textsc{WitnessPointSuperSet}}
\begin{algorithmic}
\REQUIRE  $F_1,...,F_n\in\mathbb{C}[x_1,...,x_N]$.
\ENSURE  $\{f_1,..,f_N\},\{\widehat{W}_r,..,\widehat{W}_d\}$,$L$. $\{f_1,..,f_N\}$ square system, $\widehat{W}_i$ Witness point super sets corresponding to a pure i-dimensional component of $V(f_1,...,f_N)$, $L$ a set of linear polynomials defining a linear space of dimension N-d.
\vspace{0.1cm}
\STATE $f=\{f_1,...,f_N\}$ reduction of $F=\{F_1,...,F_n\}$ to a square system (cf.\cite{Sommese2},\cite{Sommes2},\cite{Sommese6});
\STATE d=dim($V(f_1,...,f_N)$), using  Gr\"{o}bner basis (cf.\cite{Decker},\cite{Greuel});
\STATE r=N-rank(f)\footnotemark, rank(f) the rank of the Jacobian matrix of the system f at a generic point;
\STATE $L=\{l_1,...,l_d\}$ a set of $d$ generic linear polynomials;
\IF{$d=r$}
\STATE compute $T_d=V(f_1,...,f_N,l_1,...,l_d)$, using triangular sets, (cf.\cite{Decker},\cite{Greuel});
\STATE set $\widehat{W}_d=\{(x_1,...,x_N)\;\mid\;(x_1,...,x_N)\in T_d,\; (x_1,...,x_N)\in V(F)\}$;
\RETURN $\{f_1,...,f_N\},\;\{\widehat{W}_d\}$, $L$ ;
\ELSE
\FOR{$i=r\; to\; d$}
\IF{$i=0$}
\STATE $\Omega_{i}(f)(x)=f$;
\ELSE
\STATE
 \begin{equation*}
     \Omega_{i}(f)(x,z_{1},...,z_{i})=:\left(
                          \begin{array}{c}
                            f_{1}(x)+\sum_{j=1}^{i}\lambda_{1j}z_{j} \\
                            . \\
                            . \\
                            f_{N}(x)+\sum_{j=1}^{i}\lambda_{Nj}z_{j} \\
                            l_{1}+z_{1} \\
                            . \\
                            . \\
                            l_{i}+z_{i} \\
                          \end{array}
                        \right)
    \end{equation*} $\lambda_{kj}\in\mathbb{C}$ generic, $k=1,...,N$, $j=1,...,i$;
\ENDIF
\ENDFOR
\FOR{$i=d\; to\; r$}
\IF{$i=d$}
\STATE compute $T_i=V(\Omega_{i}(f)(x,z_{1},...,z_{i}))$, using triangular sets, (cf. \cite{Decker},\cite{Greuel});
\STATE set $\widehat{W}_i=\{(x_1,...,x_N)\;\mid\;(x_1,...,x_N,0,...,0)\in T_i,\; (x_1,...,x_N)\in V(F)\}$;
\STATE $S_i=T_i\setminus\{(x_1,...,x_N,z_1,...,z_i)\in T_i\;\mid\;z_1=....=z_i=0\}$;
\ELSE
\STATE compute $T_i=V(\Omega_{i}(f)(x,z_{1},...,z_{i}))$, using homotopy function with $\Omega_{i+1}(f)(x,z_{1},...,z_{i})$ as start system and $S_{i+1}$ as start solution set (cf. \cite{Sommese2},\cite{Sommes2},\cite{Sommese6});
\ENDIF
\ENDFOR
\RETURN $\{f_1,...,f_N\},\;\{\widehat{W}_r,...,\widehat{W}_d\}$, $L$ ;
\ENDIF
\end{algorithmic}
\end{algorithm}
\footnotetext[1]{ $V(f_1,...,f_N)$ has no components of dimension smaller then $N-rank(f)$ (cf. \cite{Sommese6}).}
\section{Computation Witness Point Set}\label{WitSet}
The witness point super set $\widehat{W}_i$ is an union of an i-Witness point set $W_i$ and a junk point set $J_i$ (cf. \cite{Sommese2},\cite{Sommes2},\cite{Sommese6}),
\begin{equation*}
 \widehat{W}_i=W_i\cup J_i,\;where\; W_i\subset X_i\; and\; J_i\subset \cup_{j>i}X_j\;for\; i=0,1,...,d.
\end{equation*}
We use Gr\"{o}bner bases for the 0-dimensional ideal, local dimension and homotopy continuation method to remove the points of $J_i$ from $\widehat{W}_i$ as follows:
\begin{algorithm}
\caption{\textsc{WitnessPointSet}}
\begin{algorithmic}
\REQUIRE  $\{f_1,..,f_N\}\subset\mathbb{C}[x_1,..,x_N],\;\{\widehat{W}_r,..,\widehat{W}_d\}$ a list of witness point superset sets, $L=\{l_1,..,l_d\}$ a set of generic linear polynomials (Output of \texttt{Algorithm 1}).
\ENSURE  $\{f_1,...,f_N\},\;\{W_r,..,W_d\}$, $L=\{l_1,...,l_d\}$. $W_i$ a Witness point set corresponding to a pure i-dimensional component of $V(f_1,...,f_N)$.
\vspace{0.1cm}
\STATE $W_d=\widehat{W}_d$, $s_d=\sharp{W}_d$;
\FOR{$i=d-1\; to\; r$}
\STATE $W_i=\widehat{W}_i$;
\FOR{ each point $w\in W_i$}
\STATE compute\footnotemark \; $t=dim_{w}Z$ for $Z=V(f_1-f_1(w),...,f_N-f_N(w))$, using Gr\"{o}bner basis (cf. \cite{Decker},\cite{Greuel});
\IF{$t>i$}
\STATE $W_i=W_i\setminus\{w\}$;
\ENDIF
\ENDFOR
\FOR{ each point $w\in W_i$}
\IF{$i=0$}
\STATE choose $A\subset\mathbb{C}^{d\times N}$ a generic matrix and a generic $\epsilon\in\mathbb{C}^{N}$, $\|\epsilon\|$ small;
\STATE compute $S=V(\{f_1,...,f_N,A(x-w)\})$, $T=V(\{f_1,...,f_N,A(x-w-\epsilon)\})$, using triangular sets (cf. \cite{Decker},\cite{Greuel});
\IF{$\sharp{S}=\sharp{T}$}
\STATE $W_i=W_i\setminus\{w\}$;
\ENDIF
\ELSE
\FOR{$j=i+1\; to\; d$}
\STATE choose $A\subset\mathbb{C}^{j\times N}$ a generic matrix;
\IF{$j=d$}
\STATE compute $S=V(\{f_1,...,f_N,A(x-w)\})$, using triangular sets, (cf. \cite{Decker},\cite{Greuel});
\IF{$\sharp{S}=s_d$}
\STATE $W_i=W_i\setminus\{w\}$;
\ENDIF
\ELSE
\STATE compute $S=V(\{f_1,...,f_N,A(x-w)\})$, using homotopy function with start system $\{f_1,...,f_N, l_1,...,l_j\}$ and start solution $W_j$  (cf. \cite{Sommese6});
\IF{$w\in S$}
\STATE $W_i=W_i\setminus\{w\}$;
\ENDIF
\ENDIF
\ENDFOR
\ENDIF
\ENDFOR
\ENDFOR
\RETURN $\{f_1,...,f_N\},\{W_r,...,W_d\},L$ ;
\end{algorithmic}
\end{algorithm}
\footnotetext[2]{$t=dim_{w}Z$ is less or equal to the local dimension of a point $v\in V(f_1,...,f_N)$ (cf. \cite{Fischer}) with $\|w-v\|$ small.}
\section{Partition Witness Point Sets}\label{IrrWitSet}
In this section we show that the monodromy action on an algebraic variety Z and the Zero Sum Relation are sufficient to find the breakup of the k-Witness point set $W_k$ into irreducible k-Witness point sets. We present here a modified version of the algorithms described in (cf. \cite{Sommese1},\cite{Sommese5}).\\\\
Let Z be a pure k-dimensional algebraic variety in $\mathbb{C}^N$, and $Z=\cup_{i=1}^{r}Z_i$ be the irreducible decomposition of Z. Let $\pi:\mathbb{C}^{N}\longrightarrow\mathbb{C}^k$ be a generic projection and let $l\subset\mathbb{C}^k$ be a general line.\\
Set
\begin{itemize}
 \item $\mathbb{W}_l:=\pi^{-1}(l)\cap Z$ a set of r different curves in $\mathbb{C}^N$.
 \item U is a non-empty open subset of $l$ consisting of all points $x\in l$ with $\pi^{-1}(x)$ transverse to Z.
 \item $W:=\pi^{-1}(x)\cap Z$ for a generic element $x\in U$, and $V_x:=V$ a subset of W.
 \item $W_i:=\pi^{-1}(x)\cap Z_i$ for an irreducible k-dimensional component $Z_i$ of Z.
 \item $\lambda:\mathbb{C}^N\longrightarrow\mathbb{C}$ a linear function one-to-one on $W$.
\end{itemize}
Define the function $s:U\longrightarrow\mathbb{C}$ by
\begin{equation*}
 s(y)=\sum_{z\in V_y}\lambda(z).
\end{equation*}
Where $V_y$ is a subset of $\pi^{-1}(y)\cap Z$ defined by.
\begin{equation*}
 V_y:=\{z\;\mid\; z \;on\; a\; curve\; through\; a\; point\; of\; V \}.
\end{equation*}
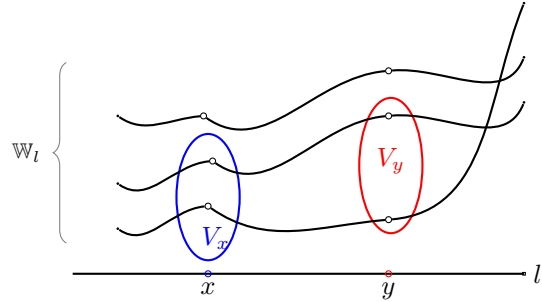
\begin{center}
\def\Piii{(-1.4,2.05) ellipse (1.5 and 0.7)} 
\tikzstyle{P_1} = [draw,red,thick]
\def\Pii{(-0.7,-2) ellipse (1.4 and 0.7)} 
\tikzstyle{P_2} = [draw,blue,thick]
\begin{tikzpicture}[scale=0.6,>=latex']
    \begin{scope}[rotate=-90]
        \path[P_1] \Piii;
        \path (-2,2.05) node[below,red]{$V_y$};
    \end{scope}
    \begin{scope}[rotate=-90]
        \path[P_2] \Pii;
        \path (0.2,-1.3) node[left,blue]{$V_x$};
    \end{scope}
\draw [thick] (-5,-1) -- (5,-1);;
\node (l) at (5,-1) [draw,inner sep=0pt] {};
\node (x) at (-2,-1) [draw,circle,blue,inner sep=0.8pt] {};
\node (y) at (2,-1) [draw,circle,red,inner sep=0.8pt] {};
\draw (x) node [below] {$x$} (l) node [right] {$l$}
   (y) node [below] {$y$};

\node (a) at (-4,0) [fill,circle,inner sep=0pt] {};
\node (x1) at (-2,0.5) [draw,circle,inner sep=0.8pt] {};
\node (x2) at (2,0.2) [draw,circle,inner sep=0.8pt] {};
\node (b) at (5,5) [fill,circle,inner sep=0pt] {};

\draw [thick] (a) to [out=-35,in=185] (x1)(x1) to [out=-35,in=185] (x2)
   (x2) to [in=248,out=5] (b);

\node (aa) at (-4,1) [fill,circle,inner sep=0pt] {};
\node (xx1) at (-1.9,1.5) [draw,circle,inner sep=0.8pt] {};
\node (xx2) at (2,2.5) [draw,circle,inner sep=0.8pt] {};
\node (bb) at (5,2.8) [fill,circle,inner sep=0pt] {};

\draw [thick] (aa) to [out=-35,in=185] (xx1)(xx1) to [out=-35,in=185] (xx2)
   (xx2) to [in=248,out=5] (bb);

\node (aaa) at (-4,2.5) [fill,circle,inner sep=0pt] {};
\node (xxx1) at (-2.1,2.5) [draw,circle,inner sep=0.8pt] {};
\node (xxx2) at (2,3.5) [draw,circle,inner sep=0.8pt] {};
\node (bbb) at (5,3.8) [fill,circle,inner sep=0pt] {};

\draw [thick] (aaa) to [out=-35,in=185] (xxx1)(xxx1) to [out=-35,in=185] (xxx2)
   (xxx2) to [in=248,out=5] (bbb);
\draw [gray,decorate,decoration={brace,amplitude=5pt},xshift=-4pt,yshift=-9pt]
   (-5,0)  -- (-5,4)
   node [black,midway,left=4pt,xshift=-2pt] {\footnotesize $\mathbb{W}_l$};
\end{tikzpicture}
\end{center}
\begin{theo} Let $l$, U, W, $W_i$ for i=1,...,r, and the functions $\lambda$, s be as above.\\
If the function s is continuous and $V\cap W_i\neq\emptyset$ for some $i\in\{1,...,r\}$, then $W_i\subseteq V$.
\end{theo}
\begin{eg} Before proving the theorem we illustrate it by an example.
\begin{itemize}
\item Let Z be the algebraic variety of dimension one in $\mathbb{C}^2$ defined by the polynomial $f(x,y)=(x^2+y^2-5)(x-2y-3)$. Let $L_1$ be the linear space of dimension one in $\mathbb{C}^2$ defined by the polynomial $l_1=x+y-3$.
\item Define a homotopy function :
 \begin{equation*}
  h(t,x(t),y(t)):=\left(
                    \begin{array}{c}
                      \alpha(t) \\
                      f(x(t),y(t)) \\
                    \end{array}
                  \right).
 \end{equation*}
 $\alpha:[0,1]\longrightarrow p^{-1}(L)$ given by
  \begin{equation*}
   \alpha(t)=(1-t)l_0+tl_1=x+y-2t-1.
  \end{equation*}
    Let $L_0$ be the 1-dimensional linear space defined by the polynomial $l_0=x+y-1$, $L_0\cap Z\neq\emptyset$. Let $G(N-k,N)$ be the Grassmanian and $R:=\{(L_{N-k},x)\in G(N-k,N)\times\overline{Z}$ $\mid x\in L_{N-k}\cap\overline{Z}\}$ be the family of the intersections $L_{N-k}\cap\overline{Z}$, $L_{N-k}\subset\mathbb{P}^N$ k-dimensional linear spaces and $\overline{Z}\subset\mathbb{P}^N$ the closure of Z. Let $p:R\longrightarrow G(N-k,N)$ be the canonical projection.
\end{itemize}
Then with conditions above $\alpha(t)$ maps a point in $L_1\cap Z$ to a point in $L_0\cap Z$ as t goes from 1 to 0.
\end{eg}
\begin{proof} (of theorem 4.1) Assume that $W_i\nsubseteq V$. Since $W_i\cap V\neq\emptyset$, then there are $a,b\in W_i$ such that a is not in V and $b\in V$. Let $a_1,...,a_r$ denote the points of the set $V\setminus\{b\}$. By (Corollary 3.5 in \cite{Sommese1}) there is a loop $\alpha$ in the fundamental group $\pi_{1}(U,\pi^{-1}(x))$ with $\alpha(0)=\alpha(1)$ which takes $a_j$ to $a_j$ for all j=1,...,r, and interchanges a and b.\\
 Since $\alpha$ is a continuous loop and $s:U\longrightarrow\mathbb{C}$ is continuous, the composition $s\circ\alpha:[0,1]\longrightarrow\mathbb{C}$ is continuous and
 \begin{equation*}
 s(\alpha(1))=s(\alpha(0))
\end{equation*}
\begin{equation*}
 \lambda(a)+\sum_{j=1}^{r}\lambda(a_j)=\lambda(b)+\sum_{j=1}^{r}\lambda(a_j),
\end{equation*}
as t goes from 1 to 0. This implies that $\lambda(a)=\lambda(b)$. But this contradicts the fact that $\lambda$ is one-to-one on W. Thus $W_i\subseteq V$.
\end{proof}
\begin{eg}
 Let Z be a pure 1-dimensional component in $\mathbb{C}^2$ defined be the polynomial $f(x,y)=(y-x)(y-2x)(y-3x)$, and $Z=Z_1\cup Z_2\cup Z_3$ be the irreducible decomposition. Let $\pi:\mathbb{C}^2\longrightarrow\mathbb{C}$ be the projection given by $\pi(x,y)=x$, and $\lambda:\mathbb{C}^2\longrightarrow\mathbb{C}$, $\lambda(x,y)=y$.\\
 Note that the restriction of $\pi$ to Z, $\pi_Z$ is proper and generically three-to-one with degree 3 equal to the degree of Z. $\lambda$ is one-to-one on the fiber $\pi^{-1}(y)=\{(x,x),(x,2x),$\\
 $(x,3x)\}$. Let L be the linear space in the Grassmannian G(1,2) defined by the linear polynomial $l(x,y)=x+y-2$. L intersects Z in the finite set $W:=\{(1,1),(\frac{2}{3},\frac{4}{3}),(\frac{1}{2},\frac{3}{2})\}$.\\
 Let $V:=\{(1,1),(\frac{2}{3},\frac{4}{3})\}\subset W$. The function $\sum_{v\in V}\lambda(v)$ given by $\lambda(x,x)+\lambda(x,2x)=x+2x=3x$ is continuous. Then by the theorem above for an irreducible 1-Witness point set $W_1\subset V$ it contains $\{(1,1)\}$.
\end{eg}
\; \textbf{Preparation of the algorithm(IrrWitnessPointSet)}:\\
 Let $Z_k=\cup_{i=1}^{r}Z_{ki}$ be the union of the irreducible k-dimensional components of the algebraic variety Z and $L_k$ be the generic linear space in $\mathbb{C}^N$ defined by k linear equations
\begin{equation*}
 l_j=c_{j0}+c_{j1}x_1+...+c_{jN}x_N.
\end{equation*}
for j=1,..,k and i=0,1,...,N , $c_{ij}\in\mathbb{C}$.\\
We use the generic linear space $L_k$ to define the generic projection\\
 $\pi:\mathbb{C}^N\longrightarrow\mathbb{C}^{k+1}$, $\pi(x_1,...,x_N):=(z_1,...,z_k,z_{k+1})$ as follows:
\begin{equation*}
 \left(
   \begin{array}{c}
     x_1 \\
     . \\
     . \\
     . \\
     . \\
     x_N \\
   \end{array}
 \right)\mapsto\left(
           \begin{array}{c}
             z_1 \\
             . \\
             . \\
             . \\
             z_{k} \\
             z_{k+1} \\
           \end{array}
         \right):=\left(
                    \begin{array}{ccccc}
                      c_{11} & c_{12} & . & . & c_{1N} \\
                      c_{21} & c_{22} & . & . & c_{1N} \\
                      . & . & . & . & . \\
                      . & . & . & . & . \\
                      c_{k1} & c_{k2} & . & . & c_{kN} \\
                      p_{1} & p_{2} & . & . & p_{N} \\
                    \end{array}
                  \right). \left(
   \begin{array}{c}
     x_1 \\
     . \\
     . \\
     . \\
     . \\
     x_N \\
   \end{array}
 \right),
\end{equation*}
 $p_1,..., p_N\in\mathbb{C}$ randomly chosen.\\
  Set $\lambda(x_1,...,x_N):=z_{k+1}$ and $l:=z_k\subset\mathbb{C}^k$ as in the theorem above. Let $y:=c_{k0}$ vary and fix the other $c_{10},...,c_{(k-1)0}$.
\begin{itemize}
\item Let $L_{k,y}$ be the linear spaces defined by the linear equations $l_1,...,l_{k-1}$ above and $l_{k,y}:=y+c_{k1}x_1+...+c_{kN}x_N$.
\item For the subset $V_y$ of $W_y:=L_{K,y}\cap Z_k$ the function s is given as a function $s_y:Z_k\cap L_{k,y}\longrightarrow\mathbb{C}$,
    \begin{equation*}
     s_y(x_1,...,x_N):=\sum_{(x_1,...,x_N)\in V_y}\lambda(x_1,...,x_N).
    \end{equation*}
 It is convenient to test the linearity of $s_y$.
\item Since $y=-(c_{k1}x_1+...+c_{kN}x_N)$, then we can define a function $s$
    \begin{equation*}
     s:\mathbb{C}\longrightarrow\mathbb{C}, s(y):=s_y(x_1,...,x_N)=\sum_{(x_1,...,x_N)\in V_y}\lambda(x_1,...,x_N).
    \end{equation*}
\item $s_y$ is a linear function in $x=(x_1,...,x_N)$ if and only if s is linear in y.
 \item To test the linearity of $s$, we take three values of y in $\mathbb{C}$, say a, b, c.\\
  If there exist $A, B\in\mathbb{C}$ such that:
    \begin{equation}
     (s(a)=Aa+B, s(b)=Ab+B)\Longrightarrow s(c)=Ac+B.
    \end{equation}
     Then s is linear in y.\\
  Here s(a), s(b) and s(c) correspond to the subsets $V_a\subset W_a=Z_k\cap L_{k,a}$, $V_b\subset W_b=Z_k\cap L_{k,b}$ and $V_c\subset W_c=Z_k\cap L_{k,c}$ respectively with $\sharp{V_a}=\sharp{V_b}=\sharp{V_c}=m$.
\end{itemize}
So far this is the approach which can be found in \cite{Sommese1}. Now we give some modifications.
\begin{itemize}
 \item The condition (4.1) of the linearity above is equivalent to the following equation
     \begin{equation}
     s(a)(b-c)+s(b)(c-a)+s(c)(a-b)=0.
    \end{equation}
\item From \texttt{Theorem 4.1} we obtain: If $W_{kj}\cap V_a\neq\emptyset$ for some $j\in\{1,...,r\}$ and the condition (4.1) of the linearity above is true, then $W_{kj}\subseteq V_a$.
 \item Let $Z(y):=\{z=\sum_{t=1}^{N}p_tv_t\mid v=(v_1,...,v_N)\in V_y,p=(p_1,...,p_N)\in\mathbb{C}^N\}$. Then
  \begin{equation*}
     s(y)=\sum_{v\in V_y}\lambda(v)=\sum_{v\in V_y}(\sum_{t=1}^{N}p_tv_t)=\sum_{z\in Z(y)}z.
    \end{equation*}
  \item To compute the sets $V_b$ and $V_c$ we use the homotopy function as t goes from 1 to 0 using $V_a$ as start set. In particular
   \begin{equation*}
    V_b:=((1-t)L_{k,b}+tL_{k,a})\cap Z, V_c:=((1-t)L_{k,c}+tL_{k,a})\cap Z,\; as\;t\;goes\; 1\rightarrow 0.
   \end{equation*}
    Continuation of the homotopy function implies that the i-th points in the sets $V_a$, $V_b$ and $V_c$ are on the same irreducible component.
  \item Let $V_a:=\{v_1,...,v_m\}$, $V_b:=\{\overline{v}_1,...,\overline{v}_m\}$ and $V_c:=\{\hat{v}_1,...,\hat{v}_m\}$ be the sets computed by using the homotopy function above . Let $Z(a):=\{a_1,...,a_m\}$, $Z(b):=\{b_1,...,b_m\}$ and $Z(c):=\{c_1,...,c_m\}$ be the sets corresponding to the set $V_a$, $V_b$ and $V_c$ respectively and defined as Z(y) above.\\\\
       From (4.2) we get a condition equivalent to the condition (4.1) of the linearity
    \begin{equation}
     (b-c)\sum_{i=1}^{m}a_i+(c-a)\sum_{i=1}^{m}b_i+(a-b)\sum_{i=1}^{m}c_i=0.
    \end{equation}
    The condition (4.3) is called \textbf{Zero Sum Relation} (cf. \cite{Corless}) of a given subset $V_a\subseteq W$ denoted by $ZSR(V_a)$.
  \item The sets $V_a$, $V_b$ and $V_c$ have distinct points and the same cardinality m, then obviously
    \begin{equation}
     ZSR(V_a)=\sum_{a_i\in V_a}ZSR(\{a_i\}).
    \end{equation}
    where $ZSR(\{a_i\})=(b-c)a_i+(c-a)b_i+(a-b)c_i$ is defined as a Zero Sum Relation of a given point in $V_a$.
\end{itemize}
\begin{algorithm}
\caption{\textsc{IrrWitnessPointSet}}
\begin{algorithmic}
\REQUIRE  $\{f_1,...,f_N\}\subset\mathbb{C}[x_1,...,x_N],\;\{W_r,...,W_d\}$, a list of witness point sets, $L=\{l_1,...,l_d\}$ a set of generic linear polynomials (Output of \texttt{Algorithm 2}). Where $W_k=\{w_1,...,w_{m_k}\}$ are witness point sets for a pure $k$-dimensional component $Z_k$ of $Z=V(f_1,...,f_N)$, $k=r,...,d$.
\ENSURE  $\{\{W_{r1},...,W_{rt_r}\},...,\{W_{d1},...,W_{dt_d}\}\}$, $W_{kr_{k}}$ irreducible Witness point sets corresponding to a $k$-dimensional irreducible component $Z_{kr_{k}}$ of $Z_k$.
\vspace{0.2cm}
\FOR{$k=r \;to\; d$}
\STATE a:=$c_{k0}$, define $L_{ka}$ to be the linear space defined by the subset $\{l_1,...,l_k\}\subset L$.
\STATE choose $b,c\in\mathbb{C}$ generic, define $L_{kb},L_{kc}$ as above;
\STATE $W_a=W_k$, $W_b=\emptyset$, $W_c=\emptyset$, $R=\emptyset$;
\STATE choose $p_1,..., p_N\in\mathbb{C}$;
\FOR{$i=1\; to\; m_k$}
\STATE compute $\{v_{i}\}\subset Z\cap L_{k,b}$ and $\{\widehat{v}_{i}\}\subset Z\cap L_{k,c}$ using the homotopy function with $\{f_1,...,f_N, l_1,...,l_{k-1},l_{k,a}\}$ as start system and $\{w_i\}$ as start solution;
\STATE compute the Zero Sum Relation of $\{w_i\}$:
 \begin{equation*}
  r_i=(a-b)(\sum_{j=1}^{N}p_j\widehat{v}_{ij})+(b-c)(\sum_{j=1}^{N}p_tw_{ij})+(c-a)(\sum_{j=1}^{N}p_tv_{ij});
 \end{equation*}
\STATE $R=R\cup\{r_i\}$\footnotemark;
\ENDFOR
\STATE int $t_k=0$;
\WHILE{$R\neq\emptyset$}
\IF{$\sum_{t\in T}t=0$ and $T$ is a smallest subset\footnotemark of $R$}
\STATE $t_k=t_k+1$;
\STATE $W_{kt_{k}}\subset W_a$ consists of the points corresponding of the points of $T$;
\STATE $R=R\setminus T$;
\ENDIF
\ENDWHILE
\ENDFOR
\RETURN  $\{\{W_{r1},...,W_{rt_r}\},...,\{W_{d1},...,W_{dt_d}\}\}$ ;
\end{algorithmic}
\end{algorithm}
\footnotetext[3]{the $i-th$ point in R corresponds to the $i-th$ point in $W_a$;}
\footnotetext[4]{smallest subset with respect to the cardinality.}
We give an example of a pure 2-dimensional variety Z which is an union of two 2-dimensional irreducible components $Z_1\; and \;Z_2$. $Z_1$ is of degree three and $Z_2$ is of degree two. The 2-Witness point set W for Z is given as a finite subset of Z consisting of five points $\{w_1,w_2,w_3,w_4,w_5\}$. $Z_1$ should contain three points $W_1:=\{w_1,w_2,w_3\}$ and the remaining points $W_2:=\{w_4,w_5\}$ are on $Z_2$. The algorithms (cf. \cite{Sommese1},\cite{Sommese5}) use the homotopy function at least nine times to breakup W into $W_1\; and\; W_2$. We will show below that we do not need more than five times to use the homotopy function to breakup W into $W_1\; and\; W_2$.
\begin{eg}.\\
 Let Z be the algebraic variety of dimension two in $\mathbb{C}^3$ defined by the polynomial $f(x,y,z)=(x^3+z)(x^2-y)$. Let L be the linear space of dimension one in $\mathbb{C}^3$ defined by the linear equations $l_1=4x+7y+2z+6,\;l_2=5x+7y+3z+6 $. Then $ W:=L\cap Z=\{w_1,w_2,w_3,w_4,w_5\}$, where\footnotemark
 \footnotetext[5]{Note that the values of $w_i$ are approximate values. The following equalities are therefore to interpret as approximations of the points $w_i$.}
$ w_1=(1,-1.1428571429,-1),w_2=(0,-0.8571428571,0),$\\
$w_3=(-0.1428571429+i*0.9147320339,-0.8163265306-i*0.2613520097,$\\
$0.1428571429-i*0.9147320339)$,\\
 $w_4=(-1,-0.5714285714,1)$,\\
 $w_5=(-0.1428571429-i*0.9147320339,-0.8163265306+i*0.2613520097,$\\
 $0.1428571429+i*0.9147320339)$.\\\\
 We now illustrate \texttt{Algorithm3\;(IrrWitnessPointSet)}:
 \begin{itemize}
  \item Use the linear space $L_1$ to define the linear projection $\pi:\mathbb{C}^3\longrightarrow\mathbb{C}^3$ as follows
   \begin{equation*}
   \pi(x,y,z):=\left(
              \begin{array}{ccc}
                4 & 7 & 2 \\
                5& 7 & 3 \\
                1 & 2 & 3\\
              \end{array}
            \right)\left(
                     \begin{array}{c}
                       x \\
                       y \\
                       z \\
                     \end{array}
                   \right)=(4x+7y+2z,5x+7y+3z,x+2y+3z).
   \end{equation*}
  \item Define the linear space $L_{1,c}$ of dimension one in $\mathbb{C}^3$ by the linear equations $l_1=4x+7y+2z+6,\;l_c=5x+7y+3z+c$, where c is generically chosen in $\mathbb{C}$. Then
      \begin{equation*}
       \pi_{Z\cap L_{1,c}}(x,y)=(-6,-c,x+2y+3z).
      \end{equation*}
  \item Define the linear function $\lambda:\mathbb{C}^2\longrightarrow\mathbb{C}$ by $\lambda(x,y,z):=x+2y+3z$.
   \item For a=6, let  $V_1=V_a:=\{w_{11}=(1,-1.1428571429,-1)\}\subset W$, $L_{1,a}:=L$ the linear space defined by $l_1=4x+7y+2z+6,\;l_a=5x+7y+3z+6$. Then\footnotemark $Z(a)=\{\sum_{v\in V_a}\lambda(v)=w_{11}[1]+2(w_{11}[2])+3(w_{11}[3])\}=\{-4.2857142858\}$.
   \item Let b=9, $L_{1,b}$ the linear space defined by $l_1=4x+7y+2z+6,\;l_b=5x+7y+3z+9$. Compute $V_b:=(tL_{1,a}+(1-t)L_{1,b})\cap Z=\{w_{12}=(1.671699881657157,-0.4776285376163331,-4.671699881657164)\}$ as t goes from 1 to 0, using $V_a$ as a start solution. $Z(b)=\{w_{12}[1]+2(w_{12}[2])+3(w_{12}[3])\}=\{-13.2986568385470012\}$.
 \footnotetext[6]{we use the notation $w_{ij}=(w_{ij}[1],w_{ij}[2],w_{ij}[3])$ for i=1,..,5, j=1,2,3.}
   \item Let c=63, $L_{1,c}$ the linear space defined by $l_1=4x+7y+2z+6,\;l_c=5x+7y+3z+63$. Compute $V_c:=(tL_{1,a}+(1-t)L_{1,c})\cap Z=\{w_{13}=(3.935100643260828,14.30425695906836,-60.93510064326094)\}$ as t goes \\
       from 1 to 0, using $V_a$ as a start solution. $Z(c)=\{w_{13}[1]+2(w_{13}[2])+3(w_{13}[3])\}=\{-150.261687368385272\}$.
  \end{itemize}
  \item Zero Sum Relation of $V_1=\{(1,-1.1428571429,-1)\}$:\\
    \begin{equation*}
     r_1:=\sum_{a\in Z(a)}(b-c)+\sum_{b\in Z(b)}(c-a)+\sum_{c\in Z(c)}(a-b)=
    \end{equation*}
            \begin{equation*}
     =-75.8098062588232524.
    \end{equation*}
    The zero sum relation set of $V_1=\{(1,-1.1428571429,-1)\}$ is $R_1:=\{r_1=-75.8098062588232524\}$.
   \begin{itemize}
   \item Let a=6, $V_a:=\{w_{11}=(0,-0.8571428571,0)\}\subset W$, $L_{1,a}:=L$ the linear space defined by $l_1=4x+7y+2z+6,\;l_a=5x+7y+3z+6$. Then $Z(a)=\{\sum_{v\in V_a}\lambda(v)=w_{11}[1]+2(w_{11}[2])+3(w_{11}[3])\}=\{-1.7142857142\}$.
   \item Let b=9, $L_{1,b}$ the linear space defined by $l_1=4x+7y+2z+6,\;l_b=5x+7y+3z+9$. Compute $V_b:=(tL_{1,a}+(1-t)L_{1,b})\cap Z=\{w_{12}=(-0.8358499408285809+i*1.046869318849985,0.2388142688081706-i*0.2991055196714253,-2.164150059171436-i*1.046869318849981)\}$ as t goes from 1 to 0, using $V_a$ as a start solution. $Z(b)=\{w_{12}[1]+2(w_{12}[2])+3(w_{12}[3])\}=\{-6.8506715807265477-i*2.6919496770428086\}$.
   \item Let c=63, $L_{1,c}$ the linear space defined by $l_1=4x+7y+2z+6,\;l_c=5x+7y+3z+63$. Compute $V_c:=(tL_{1,a}+(1-t)L_{1,c})\cap Z=$\\
       $\{w_{13}=(-1.967550321630417+i*3.257877039491183$,\\
       $15.99072866332302-i*0.9308220112831772$,\\
       $-55.03244967836969-i*3.257877039491242); \}$ as t goes from 1 to 0, using $V_a$ as a start solution. $Z(c)=\{w_{13}[1]+2(w_{13}[2])+3(w_{13}[3])\}=\{-135.083442030093447-i*8.3773981015488974\}$.
  \end{itemize}
  \item Zero Sum Relation of $V_2=\{(0,-0.8571428571,0)\}$:\\
    \begin{equation*}
     r_2:=\sum_{a\in Z(a)}(b-c)+\sum_{b\in Z(b)}(c-a)+\sum_{c\in Z(c)}(a-b)=
    \end{equation*}
        \begin{equation*}
     =107.3334745556671221-i*128.308937286793398.
    \end{equation*}
    The zero sum relation set of $V_2=\{(0,-0.8571428571,0)\}$ is\\
     $R_2:=\{r_2=107.3334745556671221-i*128.308937286793398\}$.
  \begin{itemize}
   \item For the other points $V_3=\{w_3\},V_4=\{w_4\}\; and\; V_5=\{w_5\}$, we found the zero sum relations $R_3:=\{r_3=-9.38237104997583366+i*127.0170767088\}$,\\
       $R_4:=\{r_4=-31.5236682999307779+i*128.3089372867945956\}$ and\\
       $R_5:=\{r_5=9.382371038077068-i*127.0170767088\}$.
  \item The set of Zero Sum Relation for all points of $W$ is $R=\cup_{j=1}^{5}R_j=\{r_1,r_2,r_3,r_4,r_5\}$, where i-th point in W corresponds i-th point in R.
  \item Find the smallest subset $T$ of R with $\sum_{t\in T}t=0$, which corresponds an irreducible Witness point set of W. Then we get $T_1=\{r_3,r_5\}$, $T_2=\{r_1,r_2,r_4\}$ corresponding to the irreducible Witness point sets $W_1=\{w_3,w_5\}$, $W_2=\{w_1,w_2,w_4\}$ respectively.
 \end{itemize}
\end{eg}
\begin{rem} The points of a Witness point set are computed approximately by using the homotopy continuation method. Therefore the result of the Zero Sum Relation is only almost zero.
\end{rem}
\section{Examples and timings with \singular and \bertini}
In this section we provide examples with timings of the algorithms \\
 \texttt{WitnessPointSuperSet}, \texttt{WitnessPointSet}, and \texttt{IrrWitnessPointSet} implemented in \singular to compute the numerical decomposition of a given algebraic variety defined by a polynomial system and compare them with the results of \bertini. We did not use the parallel features of \bertini.\\\\
 We tested to versions of the implementations in \bertini using the cascade algorithm and using the regenerative cascade algorithm.
Timings are conducted by using the 32-bit version of \singular{3-1-1} (cf. \cite{Decker}) and \bertini{1.2} (cf. \cite{Bates}) on an
Intel\textregistered \ Core(TM)2 Duo CPU \;\;P8400 $@$ \;2.26 GHz 2.27 GHz, 4 GB
RAM under the Kubuntu Linux operating system.\\\\\\\\
 Let Z be the algebraic variety defined by the following polynomial system:
\begin{eg} \label{ex1} (cf. \cite{Sommes2}).\\
 \begin{equation*}
  f(x,y,z)=\left(
             \begin{array}{c}
               (y-x^2)(x^2+y^2+z^2-1)(x-\frac{1}{2}) \\
               (z-x^3)(x^2+y^2+z^2-1)(y-\frac{1}{2}) \\
               (y-x^2)(z-x^3)(x^2+y^2+z^2-1)(z-\frac{1}{2}) \\
             \end{array}
           \right)
 \end{equation*}
\end{eg}
\begin{eg} \label{ex2} (cf. \cite{Sommese6},Example 13.6.4).\\
 \begin{equation*}
  f(x,y,z)=\left(
             \begin{array}{c}
               x(y^2-x^3)(x-1) \\
               x(y^2-x^3)(y-2)(3x+y)\\
             \end{array}
           \right)
 \end{equation*}
\end{eg}
\begin{eg} \label{ex3}
 \begin{equation*}
  f(x,y,z)=\left(
             \begin{array}{c}
               (x^3+z)(x^2-y) \\
               (x^3+y)(x^2-z) \\
               (x^3+z)(x^3+y)(z^2-y) \\
             \end{array}
           \right)
 \end{equation*}
\end{eg}
\begin{eg} \label{ex4}
 \begin{equation*}
  f(x,y,z)=\left(
             \begin{array}{c}
               x(y^2-x^3)(x-1) \\
               x(3x+y)(y^2-x^3)(y-2) \\
               x(y^2-x^3)(x^2-y) \\
             \end{array}
           \right)
 \end{equation*}
\end{eg}
\begin{eg} \label{ex5}
 \begin{equation*}
  f(x,y,z)=\left(
             \begin{array}{c}
               (x-1)((x^3+z)+(x^2-y)) \\
               (x^3+z)(x^2-y) \\
               (x^3+z)(x^2-1) \\
             \end{array}
           \right)
 \end{equation*}
\end{eg}

\begin{eg} \label{ex6}
 \begin{equation*}
  f(x,y,z)=\left(
             \begin{array}{c}
               (y-x^2)(x^2+y^2+z^2-1)(x-\frac{1}{2})+x^5 \\
               (z-x^3)(x^2+y^2+z^2-1)(y-\frac{1}{2})+y^4 \\
               (y-x^2)(z-x^3)(x^2+y^2+z^2-1)(z-\frac{1}{2})+z^6 \\
             \end{array}
           \right)
 \end{equation*}
\end{eg}
\begin{eg} \label{ex7}
 \begin{equation*}
  f(x,y,z)=\left(
             \begin{array}{c}
               x(y^2-x^3)(x-1)+y^2 \\
               x(y^2-x^3)(y-2)(3x+y)+x^3\\
             \end{array}
           \right)
 \end{equation*}
\end{eg}
\begin{eg} \label{ex8}
 \begin{equation*}
  f(x,y,z)=\left(
             \begin{array}{c}
               (x^3+z)(x^2-y)+x^4 \\
               (x^3+y)(x^2-z)+y^3 \\
               (x^3+z)(x^3+y)(z^2-y)+z^5 \\
             \end{array}
           \right)
 \end{equation*}
\end{eg}
\begin{eg} \label{ex9}
 \begin{equation*}
  f(x,y)=\left(
             \begin{array}{c}
               f_1=-3568891411860300072x^5+1948764938x^4+\\
               3568891411860300072x^2y^2-1948764938xy^2 \\\\
               f_2=-5105200242937540320x^5y-1701733414312513440x^4y^2+\\
               11692589628x^5+3897529876x^4y+5105200242937540320x^2y^3+\\
               1701733414312513440xy^4-11692589628x^2y^2-3897529876xy^3 \\
             \end{array}
           \right)
 \end{equation*}
\end{eg}
\begin{eg} \label{ex10}
 \begin{equation*}
  f(x,y,z)=\left(
             \begin{array}{c}
               f_1=-356737285367005125x^5-92300457164036000x^3y+\\
               1121648050080163317x^2z+290209720279281056yz \\\\
               f_2=-356737285367005125x^5+887060318883271500x^3z+\\
               1121648050080163317x^2y-2789081819567309964yz \\\\
               f_3=-356737285367005125x^5z^2+356737285367005125x^5y+\\
               887060318883271500x^3z^3-887060318883271500x^3yz+\\
               1121648050080163317x^2z^3-1121648050080163317x^2yz-\\
               2789081819567309964z^4+2789081819567309964yz^2 \\
             \end{array}
           \right)
 \end{equation*}
\end{eg}

\begin{eg}\label{ex11}
 \begin{equation*}
 f(x,y,z)=\left(
   \begin{array}{c}
    f_1= x^5y^2+2x^3y^4+xy^6+2x^3y^2z^2+2xy^4z^2+xy^2z^4-x^4y^2\\
    -2x^2y^4-y^6-x^5z-2x^3y^2z-xy^4z-2x^2y^2z^2-2y^4z^2-\\
    2x^3z^3-2xy^2z^3-y^2z^4-xz^5-3x^3y^2-3xy^4+x^4z+\\
    2x^2y^2z+y^4z-3xy^2z^2+2x^2z^3+2y^2z^3+z^5+3x^2y^2+\\
    3y^4+3x^3z+3xy^2z+3y^2z^2+3xz^3+2xy^2-3x^2z-3y^2z-\\
    3z^3-2y^2-2xz+2z\\\\
     f_2=x^6y+2x^4y^3+x^2y^5+2x^4yz^2+2x^2y^3z^2+x^2yz^4-\\
     5x^6-10x^4y^2-5x^2y^4-x^4yz-2x^2y^3z-y^5z-10x^4z^2-\\
     10x^2y^2z^2-2x^2yz^3-2y^3z^3-5x^2z^4-yz^5-3x^4y-\\
     3x^2y^3+5x^4z+10x^2y^2z+5y^4z-3x^2yz^2+10x^2z^3+\\
     10y^2z^3+5z^5+15x^4+15x^2y^2+3x^2yz+3y^3z+15x^2z^2\\
     +3yz^3+2x^2y-15x^2z-15y^2z-15z^3-10x^2-2yz+10z\\\\
     f_3=x^6y^2z+2x^4y^4z+x^2y^6z+2x^4y^2z^3+2x^2y^4z^3+\\
     x^2y^2z^5-7x^6y^2-14x^4y^4-7x^2y^6-x^6z^2-17x^4y^2z^2-\\
     17x^2y^4z^2-y^6z^2-2x^4z^4-11x^2y^2z^4-2y^4z^4-x^2z^6-\\
     y^2z^6+7x^6z+18x^4y^2z+18x^2y^4z+7y^6z+15x^4z^3+\\
     27x^2y^2z^3+15y^4z^3+9x^2z^5+9y^2z^5+z^7+21x^4y^2+\\
     21x^2y^4-4x^4z^2+13x^2y^2z^2-4y^4z^2-11x^2z^4-11y^2z^4-\\
     7z^6-21x^4z-40x^2y^2z-21y^4z-24x^2z^3-24y^2z^3-3z^5-\\
     14x^2y^2+19x^2z^2+19y^2z^2+21z^4+14x^2z+14y^2z+\\
     2z^3-14z^2\\
   \end{array}
 \right)
 \end{equation*}
\end{eg}
\begin{eg}\label{ex12}
\begin{equation*}
 f(x_1,x_2,x_3,x_4,x_5)=\left(
                          \begin{array}{c}
                            x_{5}^2+x_1+x_2+x_3+x_4-x_5-4 \\
                            x_{4}^2+x_1+x_2+x_3-x_4+x_5-4 \\
                            x_{3}^2+x_1+x_2-x_3+x_4+x_5-4 \\
                            x_{2}^2+x_1-x_2+x_3+x_4+x_5-4 \\
                            x_{1}^2-x_1+x_2+x_3+x_4+x_5-4 \\
                          \end{array}
                        \right)
\end{equation*}
\end{eg}
\begin{eg}\label{ex13}
 \begin{equation*}
  f(a,b,c,d,e,f,g)=\left(
                     \begin{array}{c}
                       a^2+2de+2cf+2bg+a\\
                       2ab+e^2+2df+2cg+b \\
                       b^2+2ac+2ef+2dg+c \\
                       2bc+2ad+f^2+2eg+d \\
                       c^2+2bd+2ae+2fg+e \\
                       2cd+2be+2af+g^2+f \\
                       d^2+2ce+2bf+2ag+g \\
                     \end{array}
                   \right)
 \end{equation*}
\end{eg}
\begin{eg}\label{ex14}. cyclic 4-roots problem.(cf.\cite{Bj1},\cite{Bj2}).
\end{eg}
\begin{eg}\label{ex15}. cyclic 5-roots problem.(cf.\cite{Bj1},\cite{Bj2}).
\end{eg}
\begin{eg}\label{ex16}. cyclic 6-roots problem.(cf.\cite{Bj1},\cite{Bj2}).
\end{eg}
\begin{eg}\label{ex17}. cyclic 7-roots problem.(cf.\cite{Bj1},\cite{Bj2}).
\end{eg}
\begin{eg}\label{ex18}. cyclic 8-roots problem.(cf.\cite{Bj1},\cite{Bj2}).
\end{eg}
\begin{eg}\label{ex19}.\\\\
 $f(x_{11},x_{12},x_{13},x_{14},x_{15},x_{21},x_{22},x_{23},x_{24},x_{25},x_{31},x_{32},x_{33},x_{34},x_{35})=$
 \begin{equation*}
  =\left(
       \begin{array}{c}
       -x_{12}x_{21}+x_{11}x_{22}\\
       -x_{13}x_{22}+x_{12}x_{23}\\
       -x_{14}x_{23}+x_{13}x_{24}\\
       -x_{15}x_{24}+x_{14}x_{25}\\
       -x_{22}x_{31}+x_{21}x_{32}\\
       -x_{23}x_{32}+x_{22}x_{33}\\
       -x_{24}x_{33}+x_{23}x_{34}\\
       -x_{25}x_{34}+x_{24}x_{35}\\
       \end{array}
  \right)
 \end{equation*}
\end{eg}
Table \ref{tabNID1} summarizes the results of the timings to compute the numerical decomposition\footnotemark.
\begin{table}[hbt]
\begin{center}
\begin{tabular}{|r|r|r|r|}
\hline
Example & \bertini & \bertini(re)& \singular\\
\hline \hline
\ref{ex1} & \hspace{1.2cm} 134.45s &  \hspace{1.2cm} 39s& \hspace{1.2cm} 36.07\\ \hline
\ref{ex2} & 3.08s & 2.5s& 1.49s \\ \hline
\ref{ex3} & 1min 21.28s & 27.4s& 4.02s \\ \hline
\ref{ex4} & 18.56s & 2.7s& 1.77s \\ \hline
\ref{ex5} & 15.36s & 8.6s& 1.29s \\ \hline
\ref{ex6} &  4min 13s &  15min 2s& 2min 27s \\ \hline
\ref{ex7} & 1.83s & 1.6s& 0.39s \\ \hline
\ref{ex8} & 3min 29s & 10min 43s& 1.69s \\ \hline
\ref{ex9} & 16s & 7s& 2s \\ \hline
\ref{ex10} & 2min 57s & 28s& 2min 35s  \\ \hline
\ref{ex11} & 44min 56s &  2min 37s& 4min 3s \\ \hline
\ref{ex12} & 4.73s & 6s& 0.37s \\ \hline
\ref{ex13} & 5.84s & 8s& 1s \\ \hline
\ref{ex14} & 1.43s & 4.3s& 0.79s \\ \hline
\ref{ex15} & 3.54s & 10s& 0.57s \\ \hline
\ref{ex16} & 3min 23.26s & 2min 29s& 1.43s \\ \hline
\ref{ex17} & 2h 11min 57s & 32min 17s& stopped after 5h \\ \hline
 \ref{ex18} & 19h48min 17s & 6h45min2s & stopped after 50h\\ \hline
  \ref{ex19} & 1min 57s & 51s& stopped after 3h  \\ \hline
\end{tabular}
\end{center}
\hspace{15mm}
\caption{Total running times for the computing a numerical decomposition of the examples above }
\label{tabNID1}
\end{table}
\footnotetext[7]{(re) means using the regenerative cascade algorithm instead of the cascade algorithm}
\begin{rem}
The timings show that for an increasing number of variables the original method of (cf.\cite{Sommese1},\cite{Sommese2},\cite{Sommes2},\cite{Sommese5},\cite{Sommese6}) becomes more efficient. One reason is that the computation of triangular sets which is used in \singular for solving polynomial systems is expensive in this case. Therefore the \texttt{Algorithm1}, \texttt{Algorithm2} become slow in this situation. This is not true for \texttt{Algorithm3}.\\
 Replacing the solving of polynomial system using triangular sets by homotopy function methods but keeping the computation of the dimension and starting in this dimension is more efficient in a case of a large number of variables.
\end{rem}

\end{document}